\documentclass[lettersize,journal]{IEEEtran}
\usepackage{amsmath,amsfonts,amssymb}
\usepackage{algorithmic}
\usepackage{algorithm}
\usepackage{array}
\usepackage{textcomp}
\usepackage{stfloats}
\usepackage{url}
\usepackage{verbatim}
\usepackage{graphicx}
\usepackage{cite}
\usepackage{mathtools}

\DeclareMathOperator{\diag}{diag}
\DeclareMathOperator{\argmax}{argmax}

\usepackage[caption=false]{subfig}

\usepackage{orcidlink}

\usepackage{arydshln}

\begin{document}

\title{Neural Network-Based Bandit: A Medium Access Control for the IIoT Alarm Scenario}

\author{Prasoon Raghuwanshi\orcidlink{0000-0002-9629-9742},~\IEEEmembership{(Student Member, ~IEEE)},
Onel Luis Alcaraz López\orcidlink{0000-0003-1838-5183},~\IEEEmembership{(Senior Member, ~IEEE)},
Neelesh B. Mehta\orcidlink{0000-0002-3614-049X},~\IEEEmembership{(Fellow, ~IEEE)},
Hirley Alves\orcidlink{0000-0002-8689-5313},~\IEEEmembership{(Member, ~IEEE)},
Matti Latva-aho\orcidlink{0000-0002-6261-0969},~\IEEEmembership{(Fellow, ~IEEE)}
\thanks{Prasoon Raghuwanshi, Onel Luis Alcaraz López, Hirley Alves, and Matti Latva-aho are with the Centre for Wireless Communications, University of Oulu, $90570$, Oulu, Finland (e-mail: Prasoon.Raghuwanshi@oulu.fi; Onel.AlcarazLopez@oulu.fi; Hirley.Alves@oulu.fi; Matti.Latva-aho@oulu.fi).}
\thanks{Neelesh B. Mehta is with the Department of Electrical Communication Engineering, Indian Institute of Science (IISc), Bangalore, 560012 India (e-mail: nbmehta@iisc.ac.in)}
\thanks{This research has been supported by the Research Council of Finland (former Academy of Finland, Grant 346208 (6G Flagship Programme) and Grant 362782), the INDIFICORE project (Grant 24650101111), the Finnish Foundation for Technology Promotion, the J. C. Bose Fellowship (Grant JCB/2023/000028), and the Riitta ja Jorma J. Takasen säätiö (Grant 20240358).}}

\maketitle

\begin{abstract}
Efficient Random Access (RA) is critical for enabling reliable communication in Industrial Internet of Things (IIoT) networks.
Herein, we propose a deep reinforcement learning-based distributed RA scheme, entitled Neural Network-Based Bandit (NNBB), for the IIoT alarm scenario.
In such a scenario, devices may detect a common critical event, and the goal is to ensure the alarm information is delivered successfully from at least one device.
The proposed NNBB scheme is implemented at each device, where it trains itself online and establishes implicit inter-device coordination to achieve the common goal.
We devise a procedure for acquiring a valuable context for NNBB, which then uses a deep neural network to process this context and let devices determine their action.
Each possible transmission pattern, i.e., transmit channel(s) allocation, constitutes a feasible action.
Our simulation results show that as the number of devices in the network increases, so does the performance gain of the NNBB compared to the Multi-Armed Bandit (MAB) RA benchmark.
For instance, NNBB experiences a 7\% success rate drop when there are four channels and the number of devices increases from 10 to 60, while MAB faces a 25\% drop.
\end{abstract}


\section{INTRODUCTION}

\IEEEPARstart{T}{he} Industrial Internet of Things (IIoT) incorporates intelligence and autonomy into manufacturing through dependable wireless connectivity, production line insights, smart alerts, and predictive maintenance via data acquisition and processing from numerous devices.
Notably, the number of IIoT devices is expected to reach around $9.4$ billion by $2030$, representing $37 \%$ of the total population of devices \cite{EUMarket} and leading to massive IIoT networks.

Massive IIoT networks face several critical challenges, including wireless channel overload, packet loss caused by high latency and deadline violation, and increased energy consumption at devices, particularly due to repeated channel access requests.
Collectively, these challenges constitute the so-called massive access problem \cite{10.1007/978-3-031-09357-9_10, 10412105}.
In general, channel access solutions addressing this problem in IIoT deployments must adapt to sporadic and correlated traffic \cite{10443958, 10323296}, require low signaling \cite{chandak2022learning, 10323296}, support decentralized decision-making \cite{10366306, 10323296}, and be scalable \cite{10366306, 10323296}.
Aloha-based Random Access (RA) schemes, such as pure Aloha and slotted Aloha, offer decentralized decision-making and are well-suited for sporadic and correlated traffic, with the added benefit of low signaling requirements.
However, they suffer from limited scalability, especially in scenarios where IIoT devices become randomly active by sensing an alarm event, i.e., leading to highly heterogeneous activation profiles over time.
This makes Aloha-based RA schemes impractical for massive IIoT networks, which together with the boom of tinyML for \textit{on-device} intelligence \cite{10285066, kallimani2024tinyml}, has motivated the research on learning-based RA for massive IIoT networks.
Learning-based RA schemes can effectively adapt to the spatial/temporal correlated nature of IIoT traffic, and even reduce the need for explicit signaling, as in \cite{chandak2022learning}.
Additionally, learning-based RA schemes can enable a decentralized modus operandi, which enhances the network's ability to tolerate faults.
In this decentralized modus operandi, each device independently learns its RA strategy, promoting scalability.

\subsection{STATE-OF-ART ON LEARNING-BASED RA}
Several learning-based RA schemes have already been proposed in the literature.
As shown in Table \ref{existing_works}, Reinforcement Learning (RL) \cite{deng2022story, deng2022reinforcement, 9415281, rech2021coordinated}, Deep RL (DRL) \cite{jiang2019distributed, 10366306, electronics12234845}, and Deep Q-Network (DQN) \cite{8665952, 9771953, s21093210, 8254101} are the most popular methods adopted by the existing learning-based RA works.

Authors in \cite{deng2022story} proposed an RL-based distributed RA scheme for Delay-Constrained (DC) communications (RLRA-DC). It establishes cooperation among devices to increase the system throughput beyond the maximum system throughput of $1/e$ achieved by the Aloha-based schemes.
However, RLRA-DC requires information about the total number of devices in the network.
Although there are approaches to estimate such information, e.g., \cite{deng2022story, lopez2022coordinated}, they are by no means error-free, inevitably leading to performance degradation of RLRA-DC in practice.
Authors in \cite{deng2022reinforcement} proposed a Tiny State-space R-learning random Access (TSRA) scheme to convey DC traffic. TSRA is a distributed RA scheme exploiting information regarding the urgency level of the data packets, which is measured in terms of their time to expire. The most important feature of TSRA is the independence of the size of the device state space from the number of devices in the network and the hard deadline for transmitting a packet. In spite of this, RLRA-DC achieves a higher system throughput than TSRA.

In \cite{9415281}, authors considered a multi-access point (AP) scenario, where each user runs a two-stage RL based RA algorithm for maximizing the throughput.
In the first stage, the user selects the AP at the beginning of the macro-time slot, while in the second stage, the user decides the sub-time slot for transmission.
Despite its appealing distributed nature, the proposal has the following limitations: (i) a large number of devices in the network would make the lookup table of the second RL algorithm too large, thus taking a significant time to converge and degrading the solution's optimality; and (ii) information about devices associated with an AP is required after every macro-time slot.

Authors in \cite{rech2021coordinated} proposed a coordinated RA scheme adapted to sporadic and correlated IIoT traffic.
Therein, the time slot selection task at each device is modelled as a Markov game, where devices are the agents, the number of retransmissions is the context, and the transmit time slots in a frame are the actions.
Moreover, each device runs an RL-based Linear Reward-Inaction (LRI) algorithm \cite{4082268} to learn the equilibrium points of the Markov game.
Despite its appealing distributed nature, LRI in \cite{rech2021coordinated} never converged to a pure Nash equilibrium, leading to the coordinated RA scheme being beneficial only for moderate traffic conditions. 
Meanwhile, authors in \cite{jiang2019distributed} developed a type of contention-cum-RL-based RA scheme called Inner-State-Driven random Access (ISDA).
ISDA uses a Neural Network (NN) with one hidden layer of five neurons to determine the transmission probabilities of devices in a time slot as a function of their inner states.
The main benefits of ISDA include the decentralized modus operandi and the ability to meet heterogeneous performance requirements.
However, the following issue needs to be addressed \cite{jiang2019distributed}: the policy optimization method is highly sub-optimal.

{
\setlength\arrayrulewidth{1pt}
\begin{table*}[!t]
\caption{Existing Learning-based RA Research Works \label{existing_works}}
\centering
\begin{tabular}{@{}p{0.55cm} p{1.5cm}@{} >{\raggedright\arraybackslash}p{5.5cm} >{\raggedright\arraybackslash}p{9cm}@{}}
\hline
\textbf{Ref.} & \textbf{Learning method} & \textbf{Scenario} & \textbf{Limitations of corresponding RA schemes}  \\
\hline
\cite{chandak2022learning} & MAB & Transmission of sensed alarm information to an external controller & $\bullet$ No more than three devices can be active concurrently \newline $\bullet$ Training occurs at the active devices but only one at a time, which is time-consuming and imposes tight synchronization requirements on the network \newline $\bullet$ MAB does not perform well with a large action space  \\
\hline
\cite{10366306} & DRL & Uplink transmission in a massive IoT network & $\bullet$ High complexity, as devices must possess a NN with $256$ neurons in each of its two hidden layers \newline $\bullet$ RA policy fails to converge when IoT devices lack information about other device's identification \\
\hline
\cite{deng2022story, deng2022reinforcement} & R-learning & Uplink transmission of delay-constrained traffic & $\bullet$ RA scheme in \cite{deng2022story} requires information about the total number of devices present in the network \newline $\bullet$ Inaccuracies in the estimation of the total number of devices lead to the performance of degradation of RA scheme in \cite{deng2022story} \newline $\bullet$ Even though RA scheme in \cite{deng2022reinforcement} does not require information about the total number of devices, however, its system throughput is lower than the one obtained in \cite{deng2022story} \\
\hline
\cite{9415281} & MAB & Uplink transmission in a multicell system & $\bullet$ Takes significant time to converge \newline $\bullet$ Information about devices associated with an AP is required after every macro-time slot \\
\hline
\cite{rech2021coordinated} & LRI & Uplink transmission of correlated IIoT traffic through a cellular network & $\bullet$ LRI never converges to a pure Nash equilibrium \newline $\bullet$ Performs badly for heavy traffic conditions  \\
\hline
\cite{jiang2019distributed} & Cross-entropy \cite{6796865} & Network having a diverse quality of service demands from multiple IoT applications & $\bullet$ The policy optimization method used in \cite{jiang2019distributed} is highly sub-optimal and it needs improvement\\
\hline
\cite{electronics12234845} & PPO & Uplink transmission in heterogeneous IoT network & $\bullet$ Requires three distinct NNs to be stored in the IoT node’s memory \newline $\bullet$ Enhances the overall network throughput without any individual performance guarantees. \\
\hline
\cite{8665952} & DQN & Coexisting networks, each operating with a distinct medium access protocol & $\bullet$ Demands an extensive NN with six hidden layers and $64$ neurons per each \newline $\bullet$ Throughput is lower than the one obtained in \cite{deng2022reinforcement} \newline $\bullet$ Requires a replica of the former NN model (also called target model), to be stored in the device’s memory, for training the latest NN model \\
\hline
\cite{9771953} & DQN & Uplink transmission of sporadic traffic on a single channel in a slotted system & $\bullet$ RA policy at each agent fails to converge due to a badly chosen reward function \newline $\bullet$ Requires a replica of the former NN model, to be stored at the central unit, for training the latest NN model  \\
\hline
\cite{s21093210} & DQN & Millimeter-wave network with numerous small/macro cells and traversing users & $\bullet$ Requires a replica of the former NN model, to be stored in the device’s memory, for training the latest NN model   \\
\hline
\cite{8254101} & DQN & Efficient spectrum access in a multichannel network with a simple collision model & $\bullet$ DQN at each device needs to be retrained offline every time a significant environmental change occurs \newline $\bullet$ Requires a replica of the former NN model, to be stored at the central unit, for training the latest NN model  \\
\hline
\end{tabular}
\end{table*}
}

Authors in \cite{10366306} proposed a multi-agent DRL-based RA scheme that employs a centralized training and decentralized execution approach for uplink transmission in massive IoT setups.
The learned RA policy in \cite{10366306} is scalable, can be shared among all devices, and can adapt to both correlated and uncorrelated IoT traffic.
However, the proposal has two flaws: (i) high complexity as devices must possess NN with $256$ neurons in each of its two hidden layers, and (ii) the RA policy fails to converge when IoT devices lack information about other device's identification.
To enhance the overall network throughput of a heterogeneous IoT network, authors in \cite{electronics12234845} proposed a distributed Proximal Policy Optimization (PPO)-based Multiple Access (MA).
IoT nodes adopting PPOMA perform online learning to update their access policies and dynamically coexist with other nodes.
A notable advantage of PPOMA is that it operates without demanding information about the MA protocols employed by the other coexisting IoT nodes. However, PPOMA requires three distinct NNs, namely actor, critic, and actor-target NNs, to be stored in the IoT node’s memory.
Additionally, \cite{electronics12234845} acknowledges that while PPOMA enhances the overall network throughput, it does so without providing any individual performance guarantees.

As mentioned earlier, DQN has also been adopted by many works.
In \cite{8665952}, authors introduced a DRL MA (DLMA) scheme exploiting DQN and aiming to maximize the sum throughput and the $\alpha$-fairness among the coexisting networks.
DLMA allows devices to learn how to properly use the time-spectral resources with no information about the medium access protocols used by other coexisting networks.
However, DLMA demands an extensive NN with six hidden layers and $64$ neurons per each, and still cannot outperform the R-learning-based TSRA in terms of throughput, as stated in \cite{deng2022reinforcement}. 
To learn a transmission policy that balances out throughput and fairness among users, authors in \cite{9771953} proposed a multi-agent DQN RA scheme.
In this work, fairness is impacted by the capture effect, which occurs when one user occupies the channel for an extended period, limiting other nodes' transmission opportunities and causing longer delays.
The capture effect is measured by the average packet age.
The scheme employs a parameter-sharing method from \cite{10.1007/978-3-319-71682-4_5} to train just one NN in a centralized manner and extends it to work for all devices.
Unfortunately, even though the RA policy at an agent can tune itself to different Poisson-distributed data arrival rates, it fails to converge because of a badly chosen reward function, as stated in \cite{jadoon2022collision}. 

Authors in \cite{s21093210} proposed a DQN-based algorithm for wandering users to estimate the congestion levels of APs and select them accordingly, thus reducing the delay experienced by the users during contention-based RA.
Meanwhile, the total number of successfully delivered packets was maximized in \cite{8254101} by solving a dynamic spectrum access problem in a multichannel network.
Due to the availability of many possible network states and the partial observability of those states due to zero message exchange among devices, a pre-trained DQN-based RA scheme is proposed.
The channel selection decisions are taken online and in a distributed fashion by devices.
However, a serious limitation is that the DQN at each device needs to be retrained offline every week or every time a significant environmental change occurs.
It is worth noting that the DQN-based works \cite{8665952, s21093210, 9771953, 8254101} have one limitation in common: they all require a replica of the former NN model (also called target model) for training the latest NN model.
Specifically, the approaches in \cite{8665952, s21093210} require a replica to be stored in the device’s memory, while those in \cite{9771953, 8254101} require it to be stored at the central unit.

\subsection{MOTIVATION AND CONTRIBUTION}\label{moti_contri}
As discussed above, there are numerous learning-based RA schemes in the recent literature addressing several wireless communication scenarios.
However, the IIoT alarm scenario remains much less explored.
The alarm scenario consists of IoT devices that may sense the same critical event, such as a gas leak, temperature anomaly, or radiation leak, and must promptly inform an External Controller (ExC).
Hence, they must adopt a low-overhead protocol that avoids collisions from simultaneous alarm transmissions so that the critical alarm message transmission from at least one device succeeds.

{
\setlength\arrayrulewidth{1pt}
\begin{table}[!t]
\caption{List of Symbols\label{table.2}}
\centering
\begin{tabular}{@{}l >{\raggedright\arraybackslash}p{6.8cm}@{}}
\hline
\textbf{Symbol} & \textbf{Description}  \\ 
\hline
$\boldsymbol{\varphi}$, $\boldsymbol{\hat{\varphi}}$ & Additive white Gaussian noise vector  \\
$\rho$ & Average transmit signal-to-noise ratio  \\
$\Lambda$ & Expected probability of successful transmission    \\
$\Omega$ & Loss function   \\
$\tau$ & Learning rate for MAB   \\
$\lambda$ & Mean-scaling multiplier   \\
$\gamma$ & Path loss exponent   \\
$\varepsilon$ & Probability of random action selection   \\
$\beta_0$ & Threshold value for ${\Vert \nabla_{\mathbf{w}} \Omega \Vert}_2$   \\
$\boldsymbol{\chi}$ & Vector for storing the clipped gradients   \\
$\mathcal{N}'$ & Batch of active devices  \\
$\mathcal{N}$ & Set of $N$ devices  \\
$\mathbf{w}$ & DNN weight vector  \\
$\mathbf{\mathring{A}}$ & Matrix containing all possible transmission patterns    \\
$\mathbf{\Psi}$ & Matrix containing the probability of each device selecting a transmission pattern    \\
$\mathbf{A}$ & Matrix storing the transmission patterns chosen by the active devices   \\
$\boldsymbol{\mathring{a}}$ & Transmission pattern   \\
$\boldsymbol{s}$ & Aggregated pilot signal received by the BS  \\
$E$ & Finite memory buffer  \\
$M$ & Number of available orthogonal channels  \\
$H$ & Number of hidden layers  \\
$h$ & Size of each hidden layer  \\
$B$ & Size of mini-batch for training the DNN  \\
$\vartheta$ & Complexity of NNBB in terms of the big-O notation \\
$\boldsymbol{a}_{\upsilon}$ & Action selected by the active agent $\upsilon$ \\
$\mathbf{c}_\upsilon$ & Channel coefficients between device $\upsilon$ and the BS  \\
$\boldsymbol{s}_\upsilon$ & Context received by the active device $\upsilon$  \\
$d_{\upsilon}$ & Distance between device $\upsilon$ and the alarm epicentre  \\
$r_{\upsilon}$ & Distance between device $\upsilon$ and the BS \\
$\vartheta_{lb}, \vartheta_{ub}$ & Lower and upper bound of NNBB's complexity \\
$\boldsymbol{\varrho_\upsilon}$ & Pilot sequence for an active agent $\upsilon$   \\
$\boldsymbol{\theta}_{\upsilon}$ & Weight vector of active agent $\upsilon$ in the case of MQLFA \\
$\varrho_{i,\upsilon}$ & $i^{th}$ pilot symbol of the active agent $\upsilon$  \\
$f(d_{\upsilon})$ & Activation probability function of device $\upsilon$  \\
$Q_{1}(\boldsymbol{\mathring{a}}_i)$ & Action value of $\boldsymbol{\mathring{a}}_i$ in the case of MAB   \\
${Q_{2}(\boldsymbol{s}_\upsilon, \boldsymbol{\mathring{a}}_i, \boldsymbol{\theta}_{\upsilon})}$ & Action value of $\boldsymbol{\mathring{a}}_i$, given $\boldsymbol{s}_\upsilon$ and $\boldsymbol{\theta}_{\upsilon}$, in the case of MQLFA \\
$\psi_{\upsilon}(\boldsymbol{\mathring{a}}_i)$ & Probability of device $\upsilon$ choosing $\boldsymbol{\mathring{a}}_i$  \\
$\xi(\mathcal{N}', \mathbf{A})$ & Indicator for a successful$/$unsuccessful transmission   \\
$\hat{q}(\boldsymbol{s},\boldsymbol{\mathring{a}}_i,\mathbf{w})$ & Parameterized action-value of $\boldsymbol{\mathring{a}}_i$ given $\boldsymbol{s}$ and $\mathbf{w}$   \\
$r(\boldsymbol{a}_{\upsilon})$ & Reward received by active agent $\upsilon$ for selecting $\boldsymbol{a}_{\upsilon}$   \\
$\boldsymbol{\phi}(\boldsymbol{s}_\upsilon, \boldsymbol{\mathring{a}}_i)$ & Feature vector of MQLFA \\
$\nabla_{\mathbf{w}} \Omega$ & Gradient of $\Omega$ with respect to $\mathbf{w}$  \\
\hline
\end{tabular}
\end{table}
}

In this work, we focus on the alarm scenario in an IIoT network. 
This problem has been studied earlier in \cite{chandak2022learning}, where authors proposed a Multi-Armed Bandit (MAB) based RA scheme that allows devices detecting the alarm to indirectly coordinate their alarm transmissions to the ExC.
The main limitations of the approach in \cite{chandak2022learning} are:
(i) it works under the assumption that no more than three devices can be active concurrently, while in practice it is not possible to control the number of active devices;
(ii) training occurs at the active devices but only one at a time, which is time-consuming and imposes tight synchronization requirements on the network; and
(iii) MAB does not perform well with a large action space.
Meanwhile, existing Industrial Wireless Sensor Network (IWSN) protocols such as WirelessHART, ISA$100.11$a, WIA-PA \cite{s24082554}, employ Time Division MA (TDMA).
TDMA is favored because it allows the prediction of communication latency for time-triggered packets in advance.
However, centralized approaches like TDMA struggle to deliver event-triggered packets, such as IIoT alarm message packets, in a timely manner since the generation of the event-triggered packets is unpredictable \cite{8737373}.
Moreover, in an IIoT alarm scenario, a large number of devices may become active simultaneously depending on device density, thus, the corresponding RA scheme must account for this.
The limitations of \cite{chandak2022learning} and existing IWSN protocols in these regards motivate our work, which presents a novel, efficient, online learning-based RA for the alarm scenario.
Our specific contributions are as follows:
\begin{itemize}
  \item We devise a novel procedure to acquire a useful context to assist the learning-based RA.
  The procedure starts with the active devices transmitting their respective pilot signals to their Base Station (BS).
  These signals are received as an aggregated signal at the BS, which then broadcasts it to devices.
  At each active device, this constitutes the context for the proposed RL framework.
  \item We propose a distributed RL-based RA scheme for the IIoT alarm scenario, referred to as Neural Network-Based Bandit (NNBB).
  NNBB uses a simple Deep NN (DNN) to process the context received by a device and decide the transmission channel(s) for the alarm signals such that the ExC successfully receives the signal on at least one channel.
  Unlike \cite{chandak2022learning}, herein there is no restriction on the number of simultaneously active devices.
  Moreover, NNBB does not require information about the total number of devices in an IIoT network.
  Note that relevant limitations for on-device context processing to the potential alternatives of DNN, such as attention-based and diffusion models, as in \cite{9770396, 10380515, 10409284}, advise against their use here. Specifically, the attention mechanism aggravates the convergence, complexity, and computational cost of the hosting RL algorithm \cite{HU2024128015}.
  Meanwhile, diffusion models require numerous denoising steps during inference to generate a single action selection probability vector.
  This not only increases the computational complexity of diffusion models but also leads to significant overhead \cite{10409284, 10081412}.
  \item We devise an online training procedure for NNBB in which the DNN weights are updated at every alarm event using Root Mean Square Propagation (RMSProp).
  Online training allows IIoT devices to adjust their actions almost immediately based on feedback from the ExC, which is crucial in industrial environments where delays can be costly.
  Additionally, the use of gradient clipping in NNBB's training procedure facilitates the smooth convergence of our multi-agent distributed RL system.
  \item We compare NNBB with three benchmark schemes, namely MAB-based RA, Random Selection (RS), and Myopic Q-learning with Linear Function Approximation (MQLFA)-based RA, under various network configurations.
  Our simulation results indicate that NNBB outperforms them even when the network experiences an increase in the available channels/devices or activation probability.
\end{itemize}

\begin{figure}[!t]
\centerline{\includegraphics[width=3.5in]{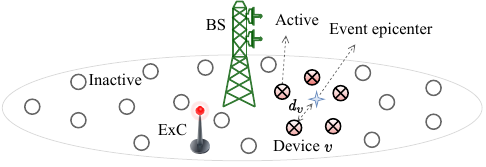}}
\caption{IIoT alarm scenario. An alarm event at a certain location triggers some nearby devices, which then become active and must convey their sensed alarm information to the ExC. Here, the BS broadcasts the context, which is the aggregated pilot signal, to the active devices.}
\label{fig.1}
\end{figure}

\subsection{ORGANIZATION}
The paper is structured as follows. Section~\ref{system_model} introduces the system model and problem formulation.
Section~\ref{NNBB_components} describes the proposed NNBB scheme and the Medium Access Control (MAC) procedure followed by an IIoT device.
Section~\ref{complexity_convergence} discusses the computational complexity and convergence of NNBB, while Section~\ref{benchmarkSchemes} discusses the benchmark schemes.
We analyze the performance of NNBB through simulations in Section~\ref{results}.
Finally, we conclude the paper in Section~\ref{conclusion} and also highlight some open research areas. 

\textit{Notation}: $\argmax(\cdot)$, $\max(\cdot)$, and $\min(\cdot)$ are the argument of the maximum function, the maximum function, and the minimum function, respectively.
$\mathcal{P}(\cdot)$ is the power set and $|\cdot|$ represents the cardinality of a set.
Superscripts $[\cdot]^T$ and $[\cdot]^H$ are transpose and conjugate transpose operations, respectively, while ${\Vert \cdot \Vert}_2$ is the Euclidean norm of a vector.
Column vectors/matrices are denoted by boldface lowercase/uppercase letters.
$\textrm{Re}(\boldsymbol{x})$ and $\textrm{Im}(\boldsymbol{x})$ return column vectors consisting of the real and imaginary components, respectively, of $\boldsymbol{x}$,
while $\textrm{abs}(\boldsymbol{x})$ returns a column vector with the absolute value of every element of $\boldsymbol{x}$.
The ${M \times M}$ identity matrix and the ${M \times 1}$ null vector are denoted as $\mathbf{I}_M$ and ${\mathbf{0}_M}$, respectively.
${\mathbb{C}^{M \times 1}}$ and ${\mathbb{R}^{M \times 1}}$ are sets of complex and real vectors of dimension $M \times 1$, respectively,
while ${\diag(\mathbf{X})}$ is a column vector comprising the diagonal elements of $\mathbf{X}$.
A circularly-symmetric complex Gaussian random vector with mean $\mathbf{\bar{y}}$ and covariance matrix $\mathbf{Z}$ is represented by ${\mathbf{y} \sim \mathcal{CN}(\mathbf{\bar{y}},\mathbf{Z})}$.
A uniform distribution between $0$ and $1$ is denoted as ${\mathcal{U}(0,1)}$.
Furthermore, ${\textrm{Pr}(X = x)}$ denotes the probability that the random variable $X$ takes on the value $x$.
Table \ref{table.2} lists the symbols used in this paper and their definitions.

\section{SYSTEM MODEL AND PROBLEM FORMULATION}      \label{system_model}

Consider the IIoT alarm scenario in Fig.~\ref{fig.1}, wherein a set $\mathcal{N}$ of wireless IIoT devices communicate with the BS and the ExC.
Specifically, IIoT devices transmit data periodically from the sensed industrial process to the BS, which consolidates it before sending it onward for processing or storage.
They are also continuously monitoring/sensing for alarm events, which they must inform to the the ExC. In this work, we focus on the timely alarm transmission process.

There are $M$ orthogonal channels for the IIoT devices to transmit alarm information concurrently.
This significantly improves the likelihood of successful reception of the information at the receiver, in contrast to the scenario in which devices have access to only a single channel.
We assume that ${|\mathcal{N}| \gg M}$ and devices are time synchronized \cite{chandak2022learning}.
At a certain time, an alarm is generated at a random location, called the epicenter.
Several examples of an alarm event are available in Section I.\ref{moti_contri}.
Devices that detect the alarm are called active devices, while the other devices are called inactive devices.
Triggered by an alarm event, a random batch of devices ${\mathcal{N}' \subseteq \mathcal{N}}$ becomes active.
Let ${f(d_{\upsilon})}$ denote the activation probability function of device ${\upsilon \in \mathcal{N}}$.
It is a decreasing function of the distance $d_{\upsilon}$ between the device $\upsilon$ and the alarm epicenter \cite{10443958, ruiz2024configuring}.

The sole purpose of the active devices is to convey the alarm information to the ExC, which is not necessarily collocated with the BS, as illustrated in Fig.~\ref{fig.1}.
The ExC responds to the received alarm information by initiating control actions to manage the industrial process.
Thus, active devices must successfully transmit the alarm information on at least one channel, irrespective of which active device makes it possible.

For transmitting the alarm information, the active device ${\upsilon \in \mathcal{N}'}$ chooses a transmission pattern denoted by ${\boldsymbol{a}_{\upsilon} = [a_{\upsilon,1}, a_{\upsilon,2}, \cdots, a_{\upsilon,M}]^T}$, where
\begin{align}
    a_{\upsilon,m} = \left\{\begin{array}{ll}
         1,& \textrm{if device} \ \upsilon \ \textrm{decides to transmit over} \\
           & \textrm{channel} \ m , \\
               0,&  \textrm{otherwise}.
               \end{array}\right.
\end{align}
We represent the transmission patterns chosen by the active devices as a matrix ${\mathbf{A}\in\{0, 1\}^{M \times | \mathcal{N}' |}}$, where each column of $\mathbf{A}$ matches with one of the vectors of the set ${\{\boldsymbol{a}_{\upsilon} | \upsilon \in \mathcal{N}'\}}$. The successful/failed alarm reception at the ExC is indicated by
\begin{align} \label{successModelequation}
    \xi(\mathcal{N}', \mathbf{A}) = \left\{\begin{array}{ll}
         \!\! 1,& \!\!\!\! \textrm{if} \ \exists m \in \{ 1,...,M \} \!\!: \!\! \displaystyle \sum_{\upsilon \in \mathcal{N}'} a_{\upsilon,m} = 1 , \\
               \!\! 0,& \!\!\!\! \textrm{otherwise},
               \end{array}\right.
\end{align}
where ${\xi(\mathcal{N}', \mathbf{A}) = 1}$ represents a successful reception of the alarm message on at least one channel.
Notice that $(\ref{successModelequation})$ does not consider decoding errors as a potential cause of transmission failures.
It only considers transmission failures resulting from medium access collision, similar to \cite{chandak2022learning, 10443958}.

\begin{figure}[!t]
\centerline{\includegraphics[width=3.5in]{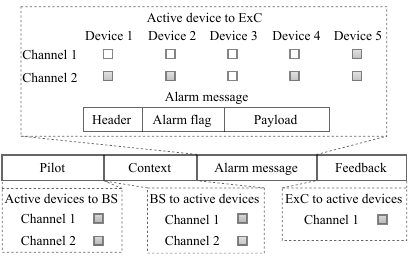}}
\caption{Illustration of the proposed protocol exemplified for the case of $M\!=$ 2 orthogonal channels.
First, each active device transmits a pilot signal to the BS.
Then, the BS broadcasts the received pilot signal to devices (context broadcast).
Finally, active devices transmit the alarm message after choosing their respective transmission patterns.
In the structure of the alarm message, the alarm flag indicates the occurrence of an event, the header carries the metadata, and the payload contains relevant information for the ExC.}
\label{alarm_message}
\end{figure}

With $M$ channels, an active device can choose from $2^M$ transmission patterns.
We compile all possible transmission patterns into a matrix ${\mathbf{\mathring{A}} = [\boldsymbol{\mathring{a}}_1, \boldsymbol{\mathring{a}}_2, \cdots, \boldsymbol{\mathring{a}}_{2^M}] \in\{0, 1\}^{M \times 2^M}}$, where each column represents a unique transmission pattern and $\boldsymbol{\mathring{a}}_i\ne \boldsymbol{\mathring{a}}_j, \forall i\ne j$.
For instance, consider the case shown in Fig.~\ref{alarm_message} where $M=2$ and $|\mathcal{N}'| = 5$.
Devices $1$, $2$, and $4$ select channel $2$, device $5$ selects both channels, and device $3$ remains silent. Then, 
\begin{align*}
\mathbf{\mathring{A}} = \left[\begin{matrix}
0 & 0 & 1 & 1\\
0 & 1 & 0 & 1
\end{matrix}\right] \textrm{and} \
\mathbf{A} = \left[\begin{matrix}
0 & 0 & 0 & 0 & 1\\
1 & 1 & 0 & 1 & 1
\end{matrix}\right]. 
\end{align*}
Therefore, $\xi(\mathcal{N}', \mathbf{A}) = 1$ since $\sum_{\upsilon \in \mathcal{N}'} a_{\upsilon,1} = 1$.

Let us define the probability of the active device $\upsilon$ choosing a transmission pattern $\boldsymbol{\mathring{a}}_i$ as $\psi_{\upsilon}(\boldsymbol{\mathring{a}}_i)$.
Let the matrix ${\mathbf{\Psi} \in[0, 1]^{|\mathcal{N}| \times 2^M}}$ store the elements of the set ${\{\psi_{\upsilon}(\boldsymbol{\mathring{a}}_i) | \forall i = 1, \cdots, 2^M\}}$ in its $\upsilon^{th}$ row.
Additionally, ${\sum_{i=1}^{2^M} \psi_{\upsilon}(\boldsymbol{\mathring{a}}_i) = 1 \textrm{,} \ \forall \upsilon \in \mathcal{N}'}$.
Now, the probability of a successful transmission can be written as \cite{chandak2022learning}
\begin{align} \label{exp_prob}
    \!\!\!\!\Lambda(\mathbf{\Psi}) = \!\!\!\!\!\!\! \sum_{\mathcal{N}' \in \mathcal{P}(\mathcal{N})} \overbrace{\prod_{\upsilon \in \mathcal{N}'} \!\! f(d_{\upsilon})}^{g_1(\mathcal{N}')} \!\!\!\!\!\!\!\!\!\!\!\! \underbrace{\sum_{\mathbf{A}\in\{0, 1\}^{M \times | \mathcal{N}' |}} \!\!\!\!\!\!\!\!\!\!\! \overbrace{\Bigl( \xi(\mathcal{N}', \mathbf{A}) \! \prod_{\upsilon \in \mathcal{N}'} \!\! \psi_{\upsilon}(\boldsymbol{a}_\upsilon) \! \Bigl)}^{g_2(\mathcal{N}', \mathbf{A})}}_{g_3(\mathcal{N}')}  .
\end{align}
Here, ${g_1(\mathcal{N}')}$ denotes the probability of devices in ${\mathcal{N}'}$ being active,
${g_2(\mathcal{N}', \mathbf{A})}$ denotes the probability of successful transmission when devices in ${\mathcal{N}'}$ chose ${\mathbf{A}}$, and
${g_3(\mathcal{N}')}$ denotes the probability of successful transmission for ${\mathcal{N}'}$.
Lastly, ${\Lambda(\mathbf{\Psi})}$ outputs the probability of successful transmission when considering all possible ${\mathcal{N}'}$.

The optimization problem for the alarm transmission task can be expressed as
\begin{subequations}\label{opti_prob}
\begin{align}
    \underset{\mathbf{\Psi}}{\textrm{maximize}} \ & \Lambda \\
    \textrm{subject to} \ & \sum_{i=1}^{2^M} \psi_{\upsilon}(\boldsymbol{\mathring{a}}_i) = 1 \textrm{,} \ \forall \upsilon \in \mathcal{N}', \\
    & \mathbf{\Psi} \in[0, 1]^{|\mathcal{N}| \times 2^M}.
\end{align}
\end{subequations}
Unfortunately, full information about ${\prod_{\upsilon \in \mathcal{N}'} f(d_{\upsilon}) , \ \forall \mathcal{N}' \in \mathcal{P}(\mathcal{N})}$ is essential for solving problem $(\ref{opti_prob})$.
This is infeasible because it requires exact knowledge about which devices become active at an alarm event.
Additionally, when ${| \mathcal{N}' | > 1}$ and/or ${|\mathcal{N}| > M}$, the problem $(\ref{opti_prob})$ becomes NP-hard, as stated in \cite{chandak2022learning}.
Therefore, we propose an online-learning-based RA scheme named NNBB for solving it.
NNBB aims to find a favorable solution $\mathbf{\Psi}^*$ without explicit information about ${\prod_{\upsilon \in \mathcal{N}'} f(d_{\upsilon})}$ and by developing an implicit coordination among devices. 
The latter assists active devices choose transmission patterns that are reliable enough to ensure a successful alarm transmission.
For this, NNBB must somehow estimate $\mathcal{N}'$ exploiting data capturing information about the identity of active devices.
We refer to this input as \textit{context}, which will allow devices with NNBB to learn its corresponding elements in $\mathbf{\Psi}^*$ autonomously.
The working principle of NNBB, including context generation and processing, is discussed next.

\section{NEURAL NETWORK-BASED BANDIT}   \label{NNBB_components}

NNBB is a DRL-based RA algorithm that is implemented in each IIoT device. The NNBB components are:
(i) an agent (device) that receives a context, (ii) actions (transmission patterns), (iii) an action value ($\psi$) corresponding to each action, and (iv) a reward (feedback signal from the ExC) for updating the action values.
The proposed protocol for the alarm transmission is illustrated in Fig.~\ref{alarm_message}.
The aforementioned NNBB components, the algorithm implementation, and its computational complexity are explained in detail below.

\begin{figure}[!t]
\centerline{\includegraphics[width=3.5in]{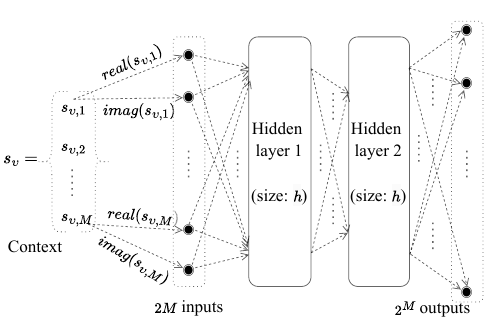}}
\caption{DNN architecture. The input to the DNN is the received context, while its outputs are 2$^M$ parameterized action values, one for each possible action.}
\label{fig.3}
\end{figure}

\subsection{CONTEXT}
An active agent (device) $\upsilon$ observes the \textit{context} $\boldsymbol{s}_\upsilon$.
For this, each active agent transmits a pilot signal consisting of $M$ pilot symbols, one on each channel.
The BS then receives the signal aggregating pilot signals from all active devices and broadcasts it to agents in the next time slot.
The signal $\boldsymbol{s}_\upsilon$ received by the active agent $\upsilon$ informs it that some devices are active in the surrounding environment and implicitly captures hidden information about their identity.
The goal of the NNBB-based RA mechanism is to make distributed access decisions based on such a context and learning experience.

Assume for simplicity that the BS, ExC, and IIoT devices are single-antenna nodes.
Let ${\boldsymbol{\varrho}_\upsilon = [\varrho_{1,\upsilon}, \dots, \varrho_{M,\upsilon}]^T \in \mathbb{C}^{M \times 1}}$ represent the pilot sequence for an active agent $\upsilon$.
Therefore, the aggregated pilot signal $\boldsymbol{s}$ received by the BS can be written as
\begin{align}
    \boldsymbol{s} = \sum_{\upsilon \in \mathcal{N}'} \sqrt{\rho} \diag(\mathbf{c}_\upsilon) \boldsymbol{\varrho}_\upsilon + \boldsymbol{\varphi} ,
\end{align}
where $\rho$ represents the average transmit signal-to-noise ratio, ${\boldsymbol{\varphi} \sim \mathcal{CN}(\mathbf{0},\mathbf{I}_M) \in \mathbb{C}^{M \times 1}}$ is the Additive White Gaussian Noise (AWGN) with normalized unit power, and ${\mathbf{c}_\upsilon \in \mathbb{C}^{M \times 1}}$ captures the $M$ channel coefficients between the device $\upsilon$ and the BS.
After receiving $\boldsymbol{s}$, the BS broadcasts it to devices.
Therefore, the context ${\boldsymbol{s}_\upsilon \in \mathbb{C}^{M \times 1}}$ received by the active device $\upsilon$ can be written as
\begin{align} \label{context_device_nu}
    \boldsymbol{s}_\upsilon = \sqrt{\rho} \diag(\mathbf{c}_\upsilon) \boldsymbol{s} + \boldsymbol{\hat{\varphi}} ,
\end{align}
where ${\boldsymbol{\hat{\varphi}} \sim \mathcal{CN}(\mathbf{0},\mathbf{I}_M) \in \mathbb{C}^{M \times 1}}$ is the AWGN with normalized power.
To simplify the analysis without any loss of generality, we assume that the transmit power of the BS is also $\rho$.

\subsection{ACTION}
After receiving $\boldsymbol{s}_\upsilon$, the active agent $\upsilon$ selects a transmission pattern $\boldsymbol{a}_{\upsilon}$, which we shall refer to as its \textit{action}, and performs a single-hop alarm transmission to the ExC.
Recall that with $M$ channels, an active device can choose from $2^M$ actions, which are collected as columns of the matrix $\mathbf{\mathring{A}}$.
The alarm message consists of a header, an alarm flag, and a payload as shown in Fig.~\ref{alarm_message}.
The alarm flag indicates the occurrence of an event, the header carries the metadata, and the payload contains relevant information about the alarm for the ExC.

\subsection{REWARD}
The active devices' cumulative target is to successfully transmit alarm information on at least one channel.
Therefore, the alarm message transmission is considered successful when the condition ${\xi(\mathcal{N}', \mathbf{A}) = 1}$ is met.
When an alarm message transmission succeeds, all active agents receive a shared Acknowledgement (ACK) from the ExC, in which case the active agent $\upsilon$ gets a reward ${r(\boldsymbol{a}_{\upsilon})=1}$.
Otherwise, ${r(\boldsymbol{a}_{\upsilon})=0}$.

\subsection{ACTION VALUE AND SELECTION}

An agent aims to maximize its total reward obtained over time.
To achieve this, it is attractive to choose those actions that have yielded a reward in the past.
This is known as exploitation. However, to discover such reward-accruing actions, an agent must try new actions, which is known as exploration.
In our work, an agent balances exploration and exploitation using the $\varepsilon$-greedy method.
This method relies on the \textit{action value} ${q(\boldsymbol{s}_\upsilon, \boldsymbol{\mathring{a}}_i), \ \forall i = 1, \cdots, 2^M}$, that represents a prediction of the expected reward for agent $\upsilon$ when it takes an action $\boldsymbol{\mathring{a}}_i$ with context $\boldsymbol{s}_\upsilon$.
In the $\varepsilon$-greedy method, an agent $\upsilon$ primarily selects a greedy action (exploitation) ${\boldsymbol{a}_{\upsilon} = \argmax_{\boldsymbol{\mathring{a}}_i \in\{0, 1\}^{M}} q(\boldsymbol{s}_\upsilon,\boldsymbol{\mathring{a}}_i)}$, with occasional exploration by randomly selecting an action $\boldsymbol{a}_{\upsilon}$ from all possible actions with a probability $\varepsilon$.
In this way, every action is eventually selected many times.
It is necessary to decrease $\varepsilon$ to gradually emphasize exploitation over exploration.
We gradually decrease $\varepsilon$ from $1$ to $0.1$ with a step size of $0.005$ after every alarm event.

Due to channel fading and noise, the number of possible contexts is infinite, leading to an agent encountering a new context every time it becomes active.
Thus, maintaining separate action values for each context is not a reasonable approach.
Rather, the agent should maintain the action values as a parameterized function and tune its parameters to better match ${q (\boldsymbol{s}_\upsilon, \boldsymbol{\mathring{a}}_i)}$ with the reward obtained after observing the feedback from the ExC.
The parameterized function, for which the DNN is employed, allows an estimation of the action value without maintaining a table of action values. 

The parameterized action-value function for the agent $\upsilon$ is expressed as ${\hat{q}(\boldsymbol{s}_\upsilon,\boldsymbol{\mathring{a}}_i,\mathbf{w})}$, where $\mathbf{w}$ is a vector of connection weights in DNN layers.
The input to the DNN is the received context, as illustrated in Fig.~\ref{fig.3}.
The DNN comprises two fully connected hidden layers and a dense output layer.
Notice that two hidden layers are capable of approximating any smooth mapping with arbitrary accuracy \cite{10.5555/1502373}.
Also, to grapple with the exploding gradients, we use gradient norm clipping \cite{gradient_clipping}.
The outputs from the DNN are $2^M$ parameterized action values $\hat{q}$.
In addition, the DNN employs: (a) ReLU non-linearity, which is a computationally efficient thresholding operation, compared to sigmoid and tanh, which helps alleviate the problem of vanishing gradients, and (b) RMSProp optimizer that maintains a moving average of the squared gradient values for each weight to adjust the update size for each weight.
Moreover, the learning rate of the DNN is a hyper-parameter in RMSProp, and we decrease it with a decay rate of $0.015$ per alarm event.

After calculating the reward ${r(\boldsymbol{a}_{\upsilon})}$, the active agent $\upsilon$ stores the tuple ${\{\boldsymbol{s}_\upsilon,\boldsymbol{a}_{\upsilon},r(\boldsymbol{a}_{\upsilon})\}}$ in its memory buffer $E$.
If $E$ is full, it removes an already stored tuple from $E$ in a First-In-First-Out (FIFO) fashion.
Indeed, there is no need to store the entire dataset in the agent's memory.
A FIFO memory buffer discards old samples to make room for new ones.
This allows an agent to be more reactive to recent changes in its environment.

The active agent $\upsilon$ trains its DNN in an online manner.
Training is performed by sampling a mini-batch of size $B$ from $E$ and providing this mini-batch as input to the DNN to update $\mathbf{w}$ by minimizing the loss function 
\begin{equation}
\Omega = \frac{1}{B} \sum_{j=1}^{B}[r_{j}(\boldsymbol{a}_\upsilon) - \hat{q}_{j}(\boldsymbol{s}_\upsilon,\boldsymbol{a}_\upsilon,\mathbf{w})]^2
\end{equation}
using RMSProp. Herein, $r_{j}$ and $\hat{q}_{j}$ represent the reward and parameterized action value, respectively, for the $j^{th}$ sample in the mini-batch.
Meanwhile, the gradient vector $\nabla_{\mathbf{w}} \Omega$ is clipped as
\begin{equation}  \label{eq5}
    \boldsymbol{\chi} \ = \ \frac{\beta_0 \ \nabla_{\mathbf{w}} \Omega}{\max({\Vert \nabla_{\mathbf{w}} \Omega \Vert}_2, \beta_0)}  ,
\end{equation}
where $\beta_0$ represents the threshold value for ${\Vert \nabla_{\mathbf{w}} \Omega \Vert}_2$ and the vector $\boldsymbol{\chi}$ stores the clipped gradients.

\begin{algorithm}[!t]
\caption{NNBB at active agent $\upsilon$}\label{Algorithm1}
\algsetup{
linenosize=\small,
linenodelimiter=.}
\begin{algorithmic}[1]
\REQUIRE Context $\boldsymbol{s}_\upsilon$, exploration rate $\varepsilon$, $\mathbf{w}$
       \STATE Evaluate ${\hat{q}(\boldsymbol{s}_\upsilon,\boldsymbol{\mathring{a}}_i,\mathbf{w}) \ \forall i = 1, \cdots, 2^M}$ using the DNN
       \STATE \label{step2} Draw $\Theta \in \mathcal{U}(0,1)$
       \IF{$\Theta > \varepsilon$} \label{step3}
          \STATE \label{step4} $\boldsymbol{a}_{\upsilon} = \argmax_{\boldsymbol{\mathring{a}}_i \in\{0, 1\}^{M}} \hat{q}(\boldsymbol{s}_\upsilon,\boldsymbol{\mathring{a}}_i,\mathbf{w})$
       \ELSE
          \STATE \label{step6} Select $\boldsymbol{a}_{\upsilon}$ randomly
       \ENDIF
       \IF{$E$ is full}
          \STATE Remove a tuple from $E$ in a FIFO fashion
       \ENDIF
       \STATE Observe the feedback from the ExC. If an ACK is received then $r(\boldsymbol{a}_{\upsilon}) = 1$, otherwise $r(\boldsymbol{a}_{\upsilon}) = 0$ 
       \STATE Store the tuple ${\{\boldsymbol{s}_\upsilon,\boldsymbol{a}_{\upsilon},r(\boldsymbol{a}_{\upsilon})\}}$ in $E$
       \STATE Sample a mini-batch of size $B$ from $E$
       \STATE Provide this mini-batch as input to the DNN for updating $\mathbf{w}$ by minimizing $\Omega$ using RMSProp
       \STATE $\varepsilon \leftarrow \max(0.1, \ \varepsilon - 0.005)$
\ENSURE Transmission pattern $\boldsymbol{a}_{\upsilon}$, $\varepsilon$, $\mathbf{w}$

\end{algorithmic}
\label{Algorithm1}
\end{algorithm}

Note that randomly sampling a mini-batch from the memory buffer breaks the temporal correlation between samples collected sequentially. 
Additionally, updating the DNN using mini-batches leads to DNN updates that are averaged over multiple samples.
This reduces the variance of DNN updates across time compared to the case where the DNN is updated with a single sample, making the training process of DNN stable.
Furthermore, using mini-batches helps prevent agents from overfitting to short-term fluctuations.
Notably, it is important to choose an appropriate mini-batch size.
A too-small mini-batch size may result in noisy DNN updates, causing oscillations in the RA policies of agents.
On the other hand, a too-large mini-batch size can make DNN updates computationally expensive and may also increase the chances of overfitting to the sampled mini-batch.

\subsection{ALGORITHM OUTLINE AND MAC ASPECTS}
The NNBB is implemented in each IIoT device detecting an alarm event. Such active devices execute the following operations:

\begin{enumerate}
\item{Transmit their respective pilot signals to the BS, which then broadcasts the received aggregated signal back to them. With this, each active agent acquires its context.}
\item{Feed their context to their DNN, which yields $2^M$ action values.}
\item{Select an action using the $\varepsilon$-greedy method and transmit the alarm message accordingly.}
\item{Monitor the feedback from the ExC and determine its reward.}
\item{Store the received context, the selected action, and the acquired reward in their memory.}
\item{Train their DNNs.}
\end{enumerate}
Algorithm \ref{Algorithm1} provides a comprehensive description of NNBB.

\begin{figure}[!t]
\centerline{\includegraphics[width=3.5in]{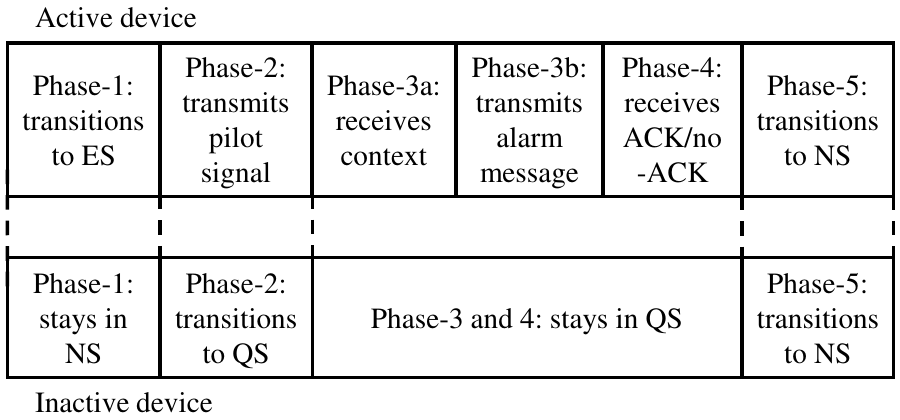}}
\caption{Illustration of the MAC procedure followed by an active and inactive device after the occurrence of an alarm event.}
\label{MAC_working}
\end{figure}

The IIoT device operation comprises three states: Normal State (NS), Emergency State (ES), and Quiet State (QS)\cite{8352005}.
In NS, the IIoT device transmits data from the sensed industrial process to the BS.
In ES, an active device transmits the alarm message to the ExC.
In QS, an inactive device halts its process data transmission to the BS.
The latter can be resumed in a future time slot, provided the device is no longer in QS.
Note that the primary duty of the BS is to collect process data from IIoT devices, while its secondary duty is to broadcast the context.

Fig.~\ref{MAC_working} illustrates the MAC procedure followed by active and inactive devices after the occurrence of an alarm event in a time slot.
During phase-1, the active devices transition to ES while the inactive devices stay in NS.
Next, in phase-2, the active devices transmit pilot signals, and the inactive devices scan available channels to determine if they are idle or busy.
If all channels are idle, the inactive devices proceed with the transmission of their process data in the subsequent phase. If all channels are busy, the inactive devices enter QS.
Following phase-2, the active devices go through phase-3a, phase-3b, and phase-4, where they receive context, transmit alarm messages, and receive ACK/no-ACK, respectively. Here, no-ACK indicates an unsuccessful alarm transmission.
Meanwhile, during phase-3, and phase-4, the inactive devices stay in QS. 
At last, in phase-5, both active and inactive devices transition to NS.

Note that our target is to improve the success rate of the first transmission attempt of the alarm message, rather than focusing on its retransmissions.
Thus, we have skipped the alarm retransmission procedure.

\section{THEORETICAL AND PRACTICAL INSIGHTS}    \label{complexity_convergence}

This section offers comprehensive discussions on the complexity, convergence, and robustness to system dynamics of the proposed NNBB.

\subsection{COMPUTATIONAL COMPLEXITY}

We quantify the computational complexity of NNBB in terms of the number of arithmetic operations executed for making a single alarm transmission decision.
The complexities of fundamental operations, such as ReLU and operations in steps~\ref{step2},~\ref{step3},~\ref{step4}, and~\ref{step6} of Algorithm~\ref{Algorithm1}, are $1,1,1, {|\mathbf{\mathring{A}}|}$, and $1$, respectively, where ${|\mathbf{\mathring{A}}| = 2^M}$.

The complexity of NNBB accounts for the complexities across three distinct phases:
action values generation phase, action selection phase, and training phase.
The complexity for the first and third phase is ${\vartheta_{1} = \sum_{i=1}^{|\boldsymbol{l}|-1} l_{i+1} (2 l_{i} + 1)}$ and ${\vartheta_{3} = B \vartheta_{1}}$, respectively, as derived in \cite{10143239}.
Here, ${\boldsymbol{l} = [l_{1}, \cdots, l_{|\boldsymbol{l}|}]^T}$ with ${l_{1} = |\boldsymbol{s}_\upsilon|}$ and ${l_{|\boldsymbol{l}|} = |\mathbf{\mathring{A}}|}$, while the remaining elements of ${\boldsymbol{l}}$ are the hidden layer sizes.
Moreover, because of steps~\ref{step3}-~\ref{step6} of Algorithm~\ref{Algorithm1}, the complexity of the second phase falls within the range ${[3, 2 + |\mathbf{\mathring{A}}|]}$.
According to Table~\ref{table.1}, we have ${\boldsymbol{l} = [M, 1, 1, 2^M]}$.
Thus, the lower bound complexity $(\vartheta_{lb})$ and upper bound complexity $(\vartheta_{ub})$ expressions for NNBB are
\begin{subequations}
\begin{align}
    \nonumber \vartheta_{lb} = & \ \vartheta_{1} + \vartheta_{3} + 3 \\
    \nonumber = & \ (B + 1) \vartheta_{1} + 3 \\
    \label{LB_NNBB} = & \ 90 \times 4^{M} + (123 + 60M)2^M + 2M + 7, \\
    \nonumber \vartheta_{ub} = & \ \vartheta_{1} + \vartheta_{3} + 2 + |\mathbf{\mathring{A}}| \\
    \label{UB_NNBB} = & \ \vartheta_{lb} + 2^M - 1.
\end{align}
\end{subequations}
By taking into account the dominant terms in $(\ref{LB_NNBB})$ and $(\ref{UB_NNBB})$, the overall complexity for NNBB is given by
\begin{align}\label{bigOComplexity}
    \vartheta & = O(4^{M}).
\end{align}
Notice that $(\ref{LB_NNBB})$, $(\ref{UB_NNBB})$, and $(\ref{bigOComplexity})$ only depend on $M$.
This means that a large value of $M$ can easily increase the complexity of NNBB beyond what resource-scarce IIoT devices can handle.
Thus, the value of $M$ in a network decides the applicability of NNBB to that network.

Consider an IIoT network with a large ${|\mathcal{N}|}$, thus, potentially large ${|\mathcal{N}'|}$ depending on the device density.
A large ${|\mathcal{N}'|}$ increases the likelihood of collisions when $M$ is small, minimizing the chances of a successful alarm transmission.
Meanwhile, a large $M$ can increase the chances of a successful alarm transmission.
However, it would also exponentially increase the complexity of NNBB, hindering NNBB's implementation on resource-scarce IIoT devices.
An appealing solution to these problems lies in:
(i) dividing the large IIoT network into several small sub-networks,
(ii) partitioning the large set of orthogonal channels into several smaller subsets, and
(iii) allocating a subset to each of these sub-networks, ensuring that the value of $M$ within each sub-network remains small.
The network division approach from \cite{10443958}, where IIoT devices with high spatial correlation are kept in different sub-networks, can be an effective method for forming sub-networks.
This way (i) NNBB implemented on an IIoT device would only take into account the value of $M$ assigned to the sub-network in which the respective IIoT device belongs, and
(ii) the chances of a successful alarm transmission would certainly increase.
By following this implementation approach, NNBB can efficiently manage alarm transmissions in large IIoT networks, without increasing its complexity.

\subsection{CONVERGENCE}

In a multi-agent distributed RL system, where multiple agents independently train their DNNs towards a common goal, ensuring convergence becomes a complex challenge.
An \textit{action-selection strategy} can significantly influence the stability and convergence of such a system by ensuring a balance between exploration and exploitation \cite{chandak2022learning}.
In this work, we are employing one of the most commonly used action-selection strategies, which is $\varepsilon$-greedy.
By carefully scheduling the decrease of $\varepsilon$, while ensuring (i) it remains high for long enough to adequately explore the action space and (ii) its minimum value is high enough to allow some level of exploration, we can prevent abrupt shifts in agents' policies, which could otherwise lead to instability.
Thus, a well-tuned $\varepsilon$-greedy strategy promotes smooth convergence in the multi-agent distributed RL system.
To satisfy aforesaid conditions, we chose to gradually decrease $\varepsilon$ from $1$ to $0.1$ with a step size of $0.005$ after every alarm event.

Meanwhile, gradient clipping plays an essential role in stabilizing the training process of the DNN, which in turn can facilitate the convergence of our multi-agent distributed RL system.
By limiting the maximum value of gradients during backpropagation, gradient clipping helps maintain a balanced learning process across all agents, preventing any single agent's updates from becoming disproportionately large, which could disrupt the collective learning process.
Additionally, large gradients can exacerbate the sensitivity to the learning rate of the DNN, potentially leading to instability if the learning rate is not finely tuned.
Gradient clipping reduces this sensitivity, making the training process of the DNN more robust to different learning rate settings.
Furthermore, gradient clipping accelerates the NN training \cite{Zhang2020Why}.
The effectiveness of gradient clipping depends on the global norm (of the gradient vector) threshold $\beta_0$.
If $\beta_0$ is too high, gradient clipping will rarely activate, causing the RA policies of agents to oscillate and struggle to converge.
Conversely, if $\beta_0$ is too low, agents may stuck with suboptimal RA policies, as the small gradient updates could stop agents from exploring the policy space.
To leverage the aforesaid benefit provided by gradient clipping, we chose to implement it in NNBB with ${\beta_0 = 5}$.

\begin{figure}[!t]
\centerline{\includegraphics[width=\columnwidth]{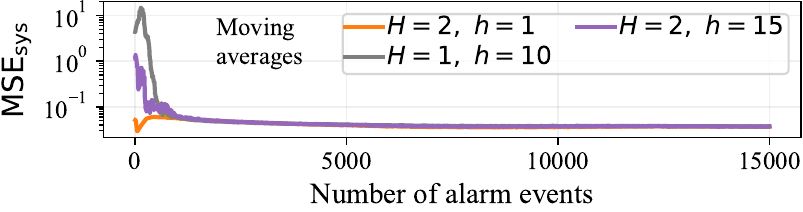}}
\caption{Variation in ${\textrm{MSE}_\textrm{sys}}$ during the training phase of NNBB for ${|\mathcal{N}|=20}$ and ${M=4}$.}
\label{convergence_N20}
\end{figure}

The mean square error of the system, at an alarm instance, is defined as
\begin{equation}
\textrm{MSE}_\textrm{sys} = \frac{1}{|\mathcal{N}'|} \sum_{\upsilon \in \mathcal{N}'} \Omega_{\upsilon} .
\end{equation}
Fig.~\ref{convergence_N20} illustrates the variation in the mean square error of the multi-agent system during the training phase of NNBB.
As can be seen in Fig.~\ref{convergence_N20}, the ${\textrm{MSE}_\textrm{sys}}$ is decreasing consistently with the progression of NNBB's training across agents, regardless of the DNN configuration.
Furthermore, ${\textrm{MSE}_\textrm{sys}}$ is becoming stable at the end for all the DNN configurations, signifying the convergence of our multi-agent distributed RL system.
Thus, Fig.~\ref{convergence_N20} proves that utilizing the right action-selection strategy and gradient clipping can facilitate the convergence of NNBB, irrespective of the system configuration.

\subsection{ROBUSTNESS TO SYSTEM DYNAMICS}
Expressions in $(\ref{context_device_nu})$, $(\ref{LB_NNBB})$, $(\ref{UB_NNBB})$, and $(\ref{bigOComplexity})$ imply that NNBB algorithm does not require information about ${|\mathcal{N}|}$.
Instead, at an alarm event, ${\mathcal{N}'}$ is the parameter of importance rather than ${|\mathcal{N}|}$.
Thus, only devices in ${\mathcal{N}'}$ are relevant to the NNBB algorithm at an alarm event.
Meanwhile, ${\mathcal{N}'}$ changes at each alarm event.
That means ${|\mathcal{N}'|}$ changes with time.
This indicates that the NNBB algorithm is robust in adapting to the varying number of devices, that are relevant at an alarm event, in the IIoT network.

Moreover, note that each agent trains its own DNN model independently in NNBB.
The DNN is designed to handle complex and time-varying inputs, meaning it could learn to recognize patterns in the input data that correlate with channel fading.
This allows agents to learn and adapt to time-varying channel fading based on their local experiences, thus, making NNBB resilient to variations in channel fading.
We have assumed Rayleigh fading for performance analysis in Section~\ref{results}, which is one of the worst types of fading due to a lack of a deterministic line-of-sight component.
Therefore, in the case of less severe (less random) fading environments, such as those with Rician or Nakagami-m $(m>1)$ fading, the performance of NNBB would certainly be better than the one obtained with Rayleigh fading.

\section{BENCHMARK SCHEMES}\label{benchmarkSchemes}


We use MAB-based RA \cite{10.5555/3312046}, RS, and MQLFA-based RA as benchmarks.
Akin to NNBB, MAB, RS, and MQLFA are implemented at each IIoT device.
In the RS scheme, actions are selected randomly at each alarm event, constituting a typical dumb approach.
Next, we comprehensively discuss the more intelligent MAB and MQLFA benchmark schemes.

\subsection{MAB-BASED RA}

MAB uses the $\varepsilon$-greedy method to select an action when an alarm event is detected.
Note that the update strategy for $\varepsilon$ in MAB-based RA is the same as the one adopted for NNBB.
The action values, ${Q_{1}(\boldsymbol{\mathring{a}}_i)}$, ${\forall i = 1, \cdots, 2^M}$, in the MAB-based RA scheme are updated as
\begin{align} \label{eq7}
  Q_{1}(\boldsymbol{\mathring{a}}_i) \leftarrow (1 - \tau) Q_{1}(\boldsymbol{\mathring{a}}_i) + r(\boldsymbol{\mathring{a}}_i) \tau .
\end{align}
Here, $\tau$ is the learning rate.
It has an initial value of $1$ and decreases with time.
In contrast with \cite{chandak2022learning}, we simultaneously train the MABs of all the active devices at each alarm event, i.e., the action values of every active device are updated.

\subsection{MQLFA-BASED RA}

MQLFA is an extension of the traditional Q-learning algorithm where each agent learns its own RA policy independently and uses the linear function approximation to compute action values as well as to handle the continuous state space.
Specifically, the action value is represented as a linear combination of features that describe the state-action pair, i.e., a weighted sum of feature values.
Thus, action values, ${Q_{2}(\boldsymbol{s}_\upsilon, \boldsymbol{\mathring{a}}_i, \boldsymbol{\theta}_{\upsilon})}$, ${\forall i = 1, \cdots, 2^M}$, in MQLFA are computed as
\begin{align} \label{MQLFA_actionValue}
  Q_{2}(\boldsymbol{s}_\upsilon, \boldsymbol{\mathring{a}}_i, \boldsymbol{\theta}_{\upsilon}) = \boldsymbol{\theta}_{\upsilon}^{T} \boldsymbol{\phi}(\boldsymbol{s}_\upsilon, \boldsymbol{\mathring{a}}_i),
\end{align}
where ${\boldsymbol{\phi}(\boldsymbol{s}_\upsilon, \boldsymbol{\mathring{a}}_i) = [\diag(\Tilde{\boldsymbol{s}}_\upsilon \Tilde{\boldsymbol{s}}_{\upsilon}^{H}), \ \boldsymbol{\mathring{a}}_i]^T \in \mathbb{R}^{2M \times 1}}$ is a feature vector, $\Tilde{\boldsymbol{s}}_\upsilon$ is the normalized $\boldsymbol{s}_\upsilon$ and is computed as
\begin{align} \label{MQLFA_nomalizedContext}
  \Tilde{\boldsymbol{s}}_\upsilon = \frac{\overbrace{\boldsymbol{s}_\upsilon - \min(\textrm{Re}(\boldsymbol{s}_\upsilon)) - \min(\textrm{Im}(\boldsymbol{s}_\upsilon))j}^{\boldsymbol{\kappa}}}{\max(\textrm{abs}(\boldsymbol{\kappa}))}.
\end{align}
Moreover, ${\boldsymbol{\theta}_{\upsilon} \in \mathbb{R}^{2M \times 1}}$ is a weight vector and each agent learns its ${\boldsymbol{\theta}_{\upsilon}}$ independently.

In MQLFA-based RA, after detecting an alarm, active agents first transmit their respective pilot signals and receive the context $\boldsymbol{s}_\upsilon$.
Then, active agents compute their respective action values, ${Q_{2}(\boldsymbol{s}_\upsilon, \boldsymbol{\mathring{a}}_i, \boldsymbol{\theta}_{\upsilon})}$, ${\forall i = 1, \cdots, 2^M}$, ${\forall \upsilon \in \mathcal{N}'}$,
and utilize them to select their respective actions, $\boldsymbol{a}_{\upsilon}$, ${\forall \upsilon \in \mathcal{N}'}$, with the help of the $\varepsilon$-greedy method.
Note that the update strategy for $\varepsilon$ in MQLFA-based RA is the same as the one adopted for NNBB.
Next, after receiving the reward ${r(\boldsymbol{a}_{\upsilon})}$, active agents update their respective ${\boldsymbol{\theta}_{\upsilon}}$, ${\forall \upsilon \in \mathcal{N}'}$, as
\begin{align} \label{MQLFA_weightUpdate}
  \boldsymbol{\theta}_{\upsilon} \leftarrow \boldsymbol{\theta}_{\upsilon} + (r(\boldsymbol{a}_{\upsilon}) - Q_{2}(\boldsymbol{s}_\upsilon, \boldsymbol{a}_{\upsilon}, \boldsymbol{\theta}_{\upsilon})) \tau \boldsymbol{\phi}(\boldsymbol{s}_\upsilon, \boldsymbol{a}_{\upsilon}) .
\end{align}
Note that MQLFA-based RA takes into account $\boldsymbol{s}_\upsilon$.
This makes the MQLFA-based RA scheme a suitable benchmark for our proposed NNBB scheme.

{
\setlength\arrayrulewidth{1pt}
\begin{table}[!t]
\caption{Parameters Used in Simulations\label{table.1}}
\centering
\begin{tabular}{p{6cm} l}
\hline
\textbf{Parameters} & \textbf{Value}  \\ 
\hline
Density of the circular region & $0.2$ devices$/$m$^2$  \\
Initial learning rate for DNN and MAB & $1.0$  \\
Total number of algorithm runs & $100$  \\
Mini-batch size $(B)$ & $2^M \times 30$  \\
Memory buffer size $(E)$ \cite{nabati2021online} & $2^M \times 100$  \\
Number of hidden layers $(H)$ & $2$   \\
Size of each hidden layer $(h)$ & $1$   \\
Global norm (of gradient vector) threshold $(\beta_0)$ \cite{nabati2021online} & $5.0$   \\
Path loss exponent $(\gamma)$ & $3.8$   \\
Mean-scaling multiplier $(\lambda)$ & $3$  \\
Initial value of ${\boldsymbol{\theta}_{\upsilon}}$ & $\mathbf{0}_{2M}$  \\
\hline
\end{tabular}
\end{table}
}

\section{NUMERICAL RESULTS}    \label{results}

We now assess the performance of the proposed NNBB RA scheme through simulations.
We consider a circular region around the BS, where the ExC and devices are uniformly distributed.
The circular region has a fixed device density of $0.2$ devices$/$m$^2$, and we vary the radius of the area to accommodate the specified number of devices.
Similar to \cite{10443958}, we use
\begin{align} \label{eq6} 
    f(d_{\upsilon}) = e^{-d_{\upsilon}/\lambda},
\end{align}
as the activation probability function, where ${\lambda >0}$ is a mean-scaling multiplier.

Pilot symbols ${\{\varrho_{i,\upsilon} | \forall i = 1, \cdots, M\}}$ are sampled from the symmetric quadrature phase shift keying constellation; i.e., chosen uniformly at random from ${\mathcal{B} = \{\frac{1+j}{\sqrt{2}}, \frac{-1-j}{\sqrt{2}}, \frac{1-j}{\sqrt{2}}, \frac{-1+j}{\sqrt{2}}\}}$.
Note that pilot symbols ${\{\varrho_{i,\upsilon} | \forall i = 1, \cdots, M\}}$ are independent and identically distributed random variables with probability mass function ${\{\textrm{Pr}(\varrho_{i,\upsilon} = x) = \frac{1}{4} | x \in \mathcal{B}, \forall i = 1, \cdots, M\}}$.
The channel undergoes quasi-static Rayleigh fading and ${\mathbf{c}_\upsilon \sim \mathcal{CN}(\mathbf{0}, \frac{1}{r_{\upsilon}^{\gamma}} \mathbf{I}_M)}$, where ${r_{\upsilon}^{\gamma}}$ represents the path loss, $\gamma$ is the path loss exponent, and $r_{\upsilon}$ is the distance between the device $\upsilon$ and the BS.
The simulation parameters and their corresponding values are listed in Table~\ref{table.1}.

We measure the performance of an RA scheme using the success rate metric, which is defined as the probability that an alarm event is successfully reported to the ExC.
We measure the success rate of an algorithm after it has converged.
To illustrate the variability in the success rate of the first transmission attempt in our simulation results, we generate an alarm at each time slot and discard the respective alarm message from the active devices after its first transmission attempt.

\begin{figure}[!t]
\centerline{\includegraphics[width=\columnwidth]{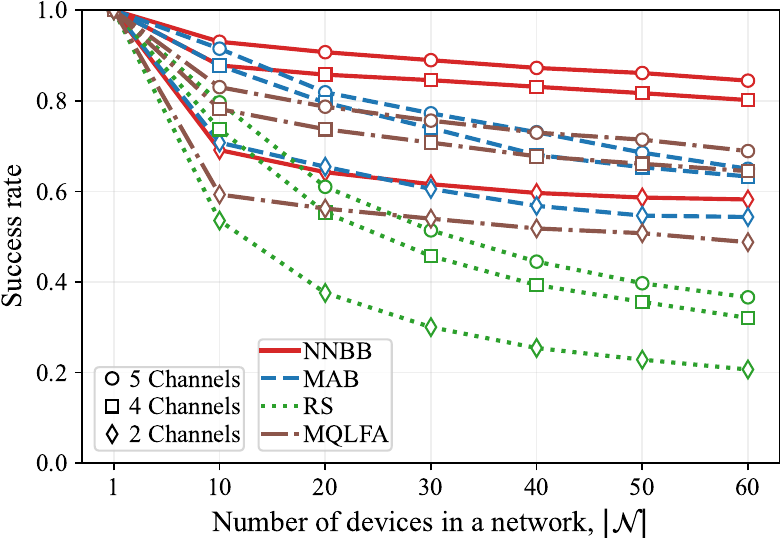}}
\caption{Comparison of the success rates of various RA algorithms for different numbers of available channels in the network.}
\label{fig.4}
\end{figure}

\subsection{IMPACT OF THE NUMBER OF DEVICES AND CHANNELS}

Fig.~\ref{fig.4} shows the impact of the number of devices on the success rate achieved by NNBB and benchmark schemes.
As the number of devices increases, the success rate deteriorates.
However, NNBB increasingly outperforms MAB, RS, and MQLFA.
Specifically, for a network of $20$ devices with $5$ channels, the success rate of NNBB is $9 \%$ higher than that of MAB.
This gain increases to $20 \%$ when there are $60$ devices.
Although MAB outperforms NNBB when only $2$ channels are shared by at most $20$ devices, the superiority is marginal.
The superiority disappears and gets reversed as the number of devices increases.
Meanwhile, NNBB outperforms MQLFA for all network configurations.
Even though both NNBB and MQLFA utilize context, NNBB processes it more effectively due to the use of the DNN.
This allows NNBB to compute relatively better action values than MQLFA, which later leads to relatively better action selection.
Moreover, the utilization of context in NNBB and MQLFA makes them more robust to an increase in the number of devices compared to MAB and RS.
This allows MQLFA to surpass MAB in terms of success rate when the number of devices goes beyond $40$ in some network configurations.
All in all, as the number of devices increases, the performance gain of NNBB over MAB, RS, and MQLFA increases.

\begin{figure}[!t]
\centerline{\includegraphics[width=\columnwidth]{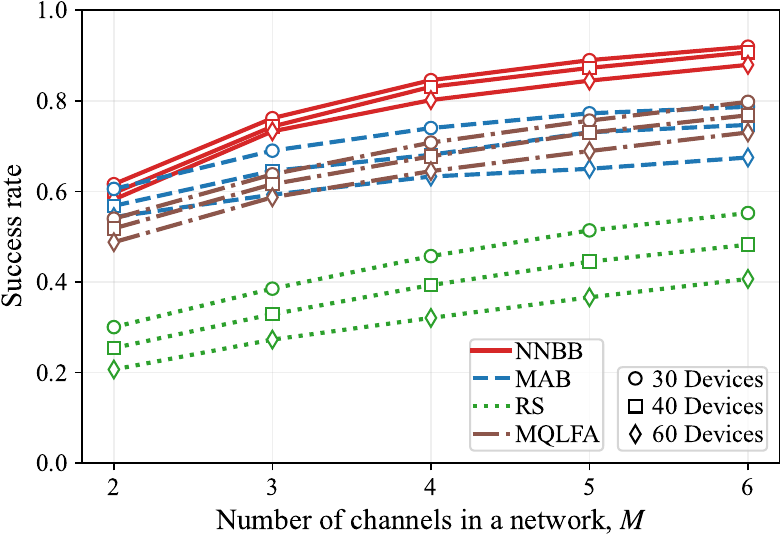}}
\caption{Comparison of the success rates of various RA algorithms for different numbers of devices in the network.}
\label{fig.5}
\end{figure}

Fig.~\ref{fig.5} plots the success rate as a function of the number of channels.
We observe that NNBB outperforms MAB, RS, and MQLFA regardless of the number of available channels for heavily loaded networks.
Notably, increasing the number of available channels from $2$ to $6$ results in a $30 \%$ increase in the success rate for NNBB in a network of $60$ devices, whereas MAB only experiences a $13 \%$ improvement.
A closer examination of the plots in Fig.~\ref{fig.5} also reveals that NNBB is less sensitive to an increase in the number of devices than MAB, RS, and MQLFA.
Moreover, the exploitation of context in NNBB and MQLFA is making them more robust than MAB and RS to an increase in the number of channels.

\begin{figure}[!t]
\centerline{\includegraphics[width=\columnwidth]{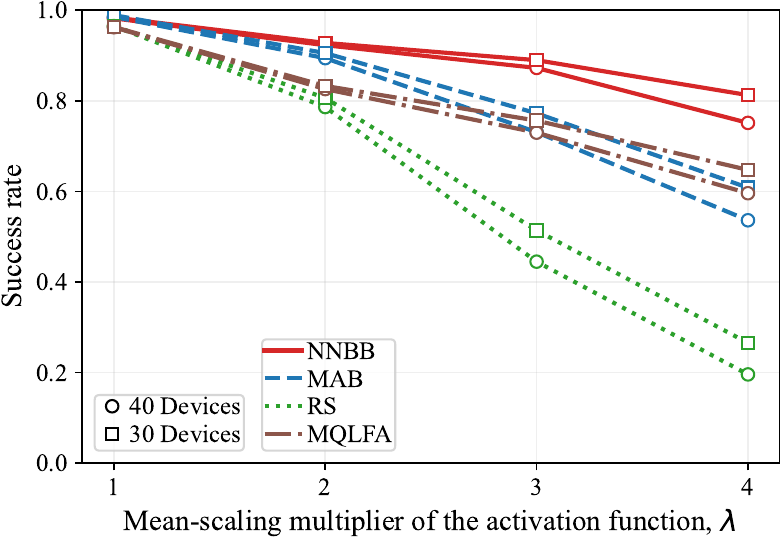}}
\caption{Comparison of the success rates of various RA algorithms as a function of the scaling multipliers $\lambda$, when the network consists of 5 channels.}
\label{fig.6}
\end{figure}

\subsection{IMPACT OF THE ACTIVATION PROBABILITY}

Fig.~\ref{fig.6} illustrates the effect of $\lambda$, which directly impacts the activation probability, on the performance of the RA schemes.
As $\lambda$ increases, the activation probability increases.
The figure indicates that NNBB and benchmark schemes exhibit similar success rates when the activation probability is small, i.e., when $\lambda=1$.
However, as the activation probability increases, the success rate of NNBB decays at a slower rate compared to MAB, RS, and MQLFA.
Even for $\lambda=3$, when the success rates of benchmark schemes fall well below $0.8$, the success rate of NNBB remains significantly higher.
Notably, for NNBB, a significant difference in the success rates between the two considered networks is visible only when $\lambda=4$, while for benchmark schemes the difference becomes apparent already when $\lambda=3$.

\subsection{IMPACT OF THE NUMBER AND SIZE OF THE HIDDEN LAYERS}

\begin{figure}[!t]
\centerline{\includegraphics[width=\columnwidth]{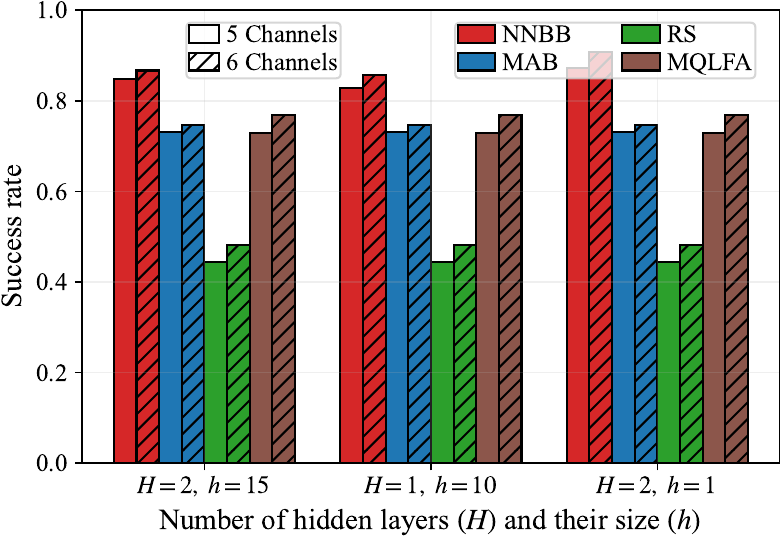}}
\caption{Comparison of the success rates of NNBB for various combinations of hidden layers, when the network consists of 40 devices.}
\label{fig.7}
\end{figure}

\begin{figure}[!t]
\captionsetup[subfigure]{labelformat=empty}
\centering
\begin{minipage}[t]{0.49\columnwidth}
\includegraphics[width=\linewidth]{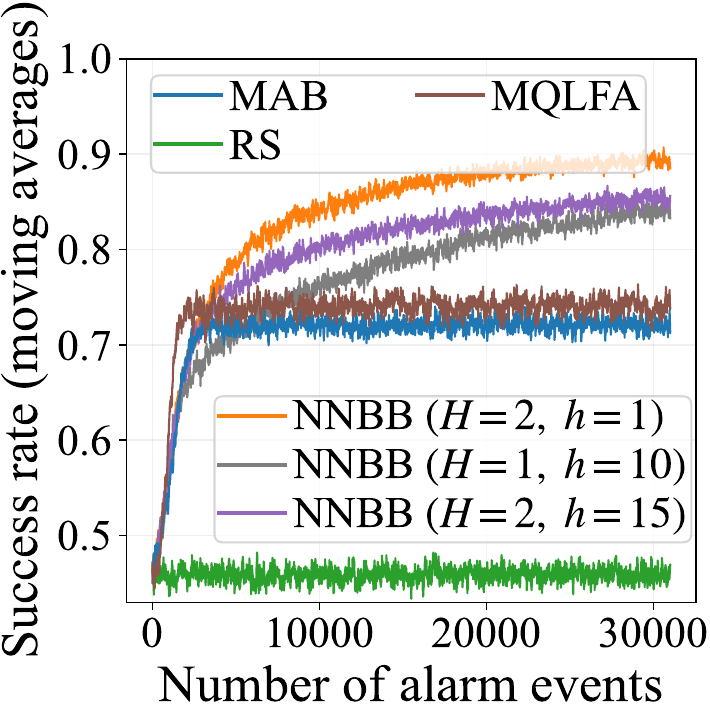}
\label{fig.8(a)}
\end{minipage}\hfill
\begin{minipage}[t]{0.49\columnwidth}
\includegraphics[width=\linewidth]{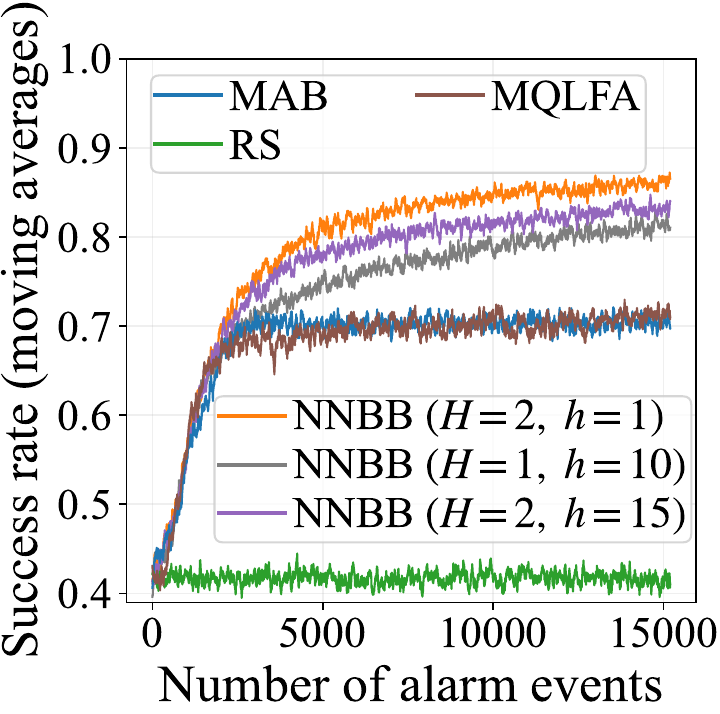}
\label{fig.8(b)}
\end{minipage}
\caption{Comparison of the performance of NNBB (during training) for various combinations of hidden layers, when the network consists of 40 devices and (a) 6 channels (left) or (b) 5 channels (right).}
\label{fig.8}
\end{figure}

Fig.~\ref{fig.7} illustrates the performance of NNBB as a function of the number of hidden layers and their size.
It is evident from the figure that NNBB consistently outperforms the benchmark approaches.
Interestingly, the hidden layer combination with the least number of neurons, i.e., ${(H, h)=(2, 1)}$, yields the highest success rate among the considered configurations.
Additionally, this configuration outperforms the other two throughout the training procedure, as indicated in Fig.~\ref{fig.8}.
It is often advantageous to reduce the hidden layer size as much as possible to maintain the NN's generalization capability.
Excess neurons can cause the network to act like a memory bank, leading to suboptimal performance when presented with input other than training samples \cite{post}.
This justifies the suboptimal performance of the configurations ${(H, h)=(1, 10)}$ and ${(H, h)=(2, 15)}$.
Reducing the hidden layer size also reduces the computational complexity of NNBB.
Furthermore, Fig.~\ref{fig.8} demonstrates that for a network of $40$ devices, NNBB should train over $15200$ and $31000$ alarm events to converge when $5$ and $6$ channels, respectively, are available.

\section{CONCLUSION}    \label{conclusion}
We proposed NNBB, a distributed DRL-based RA scheme that allows IIoT devices to develop implicit coordination to successfully convey an alarm message to an external controller.
Specifically, in the proposed scheme, upon the detection of an alarm event, every active device starts a procedure to acquire a useful context to feed a local DNN.
Then, with the help of the DNN and the $\varepsilon$-greedy method, the device selects the transmit channel(s) for the alarm message, including also the possibility of no transmission.
A reward or penalty is granted based on the success or failure of the transmission, which eventually is used to train the DNN.
Notably, the proposed DNN employs two hidden layers with just one neuron each, significantly decreasing NNBB’s computational complexity.
Simulations revealed that NNBB experiences a relatively lower drop in its success rate compared to the benchmark schemes as the number of devices in a network increases, while the success rate gained by increasing the available channels is relatively higher.

For future work, it might be interesting to substitute the adopted $\varepsilon$-greedy method, which is less effective in large action spaces and sensitive to the initial value of $\varepsilon$, by Thompson sampling-based exploration or upper confidence bound-based exploration.
Meanwhile, NNBB with the help of context develops implicit coordination among devices.
Therefore, exploiting NNBB in connected robotics and autonomous systems is an interesting future research direction.

\bibliographystyle{IEEEtran}
\bibliography{IEEEabrv,references}

\begin{IEEEbiography}[{\includegraphics[width=1in,height=1.25in,clip,keepaspectratio]{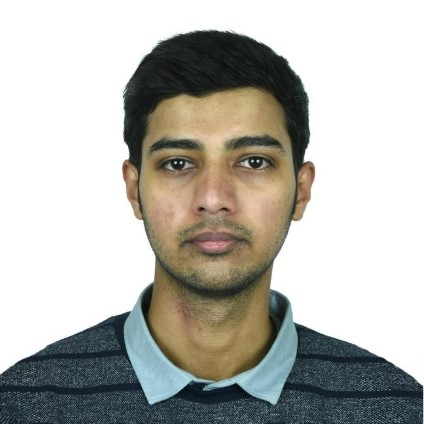}}]{Prasoon Raghuwanshi } (Student Member, IEEE) received the B.Tech and M.Tech degrees in electronics $\&$ communication engineering from the National Institute of Technology Hamirpur, Hamirpur, India, in 2020. He is currently pursuing a D.Sc. degree in communications engineering from the University of Oulu, Oulu, Finland. He joined the Centre for Wireless Communications, University of Oulu, in 2022. His research interests include random access protocols for IoT networks, machine-type communications, goal-oriented communications, deep reinforcement learning, and TinyML. \par
\end{IEEEbiography}

\begin{IEEEbiography}[{\includegraphics[width=1in,height=1.25in,clip,keepaspectratio]{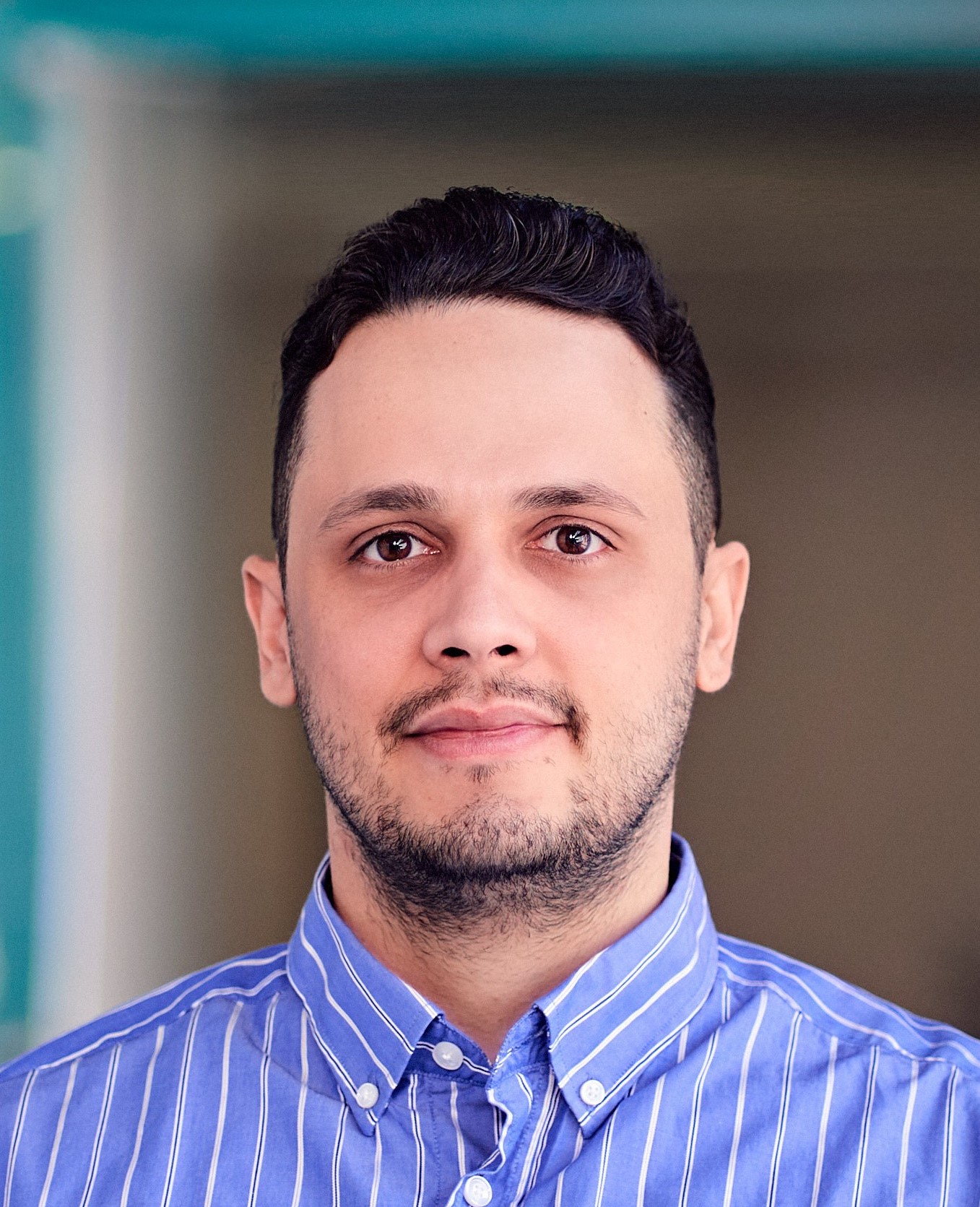}}]{Onel L. A. L\'opez }  (Senior Member, IEEE) (S'17-M'20-SM'24) received the B.Sc. (1st class honors, 2013), M.Sc. (2017), and D.Sc. (with distinction, 2020) degree in Electrical Engineering from the Central University of Las Villas (Cuba),  the Federal University of Paran\'a (Brazil), and the University of Oulu (Finland), respectively.  

He is a collaborator to the 2016 Research Award given by the Cuban Academy of Sciences, a co-recipient of the 2019 and 2023 IEEE EuCNC Best Student Paper Award, and the recipient of the 2020 Best Doctoral Thesis Award granted by Finland TEK and TFiF in 2021. He is co-author of the books entitled ``Wireless RF Energy Transfer in the Massive IoT Era: Towards Sustainable Zero-energy Networks'', Wiley, 2021, and "Ultra-Reliable Low-Latency Communications: Foundations, Enablers, System Design, and Evolution Towards 6G'', Now Publishers, 2023. He is currently an Associate Professor (tenure track) in sustainable wireless communications engineering at the Centre for Wireless Communications (CWC), Oulu, Finland. His research interests include wireless communications, signal processing, sustainable IoT, and wireless RF energy transfer.\par 
\end{IEEEbiography}

\begin{IEEEbiography}[{\includegraphics[width=1in,height=1.25in,clip,keepaspectratio]{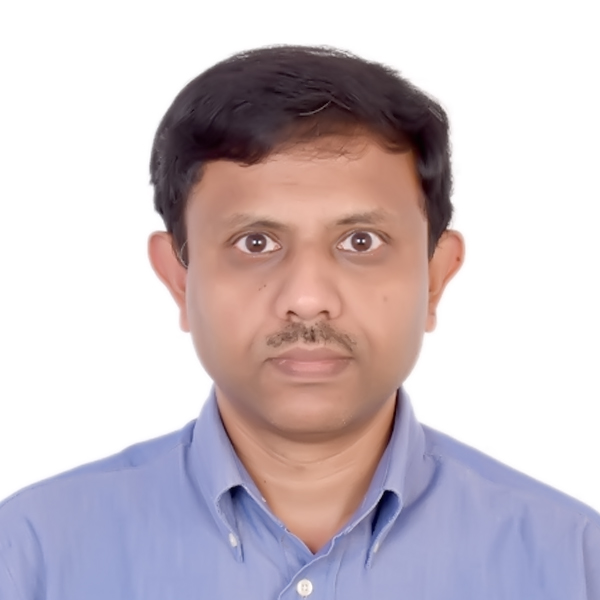}}]{Neelesh B. Mehta } (Fellow, IEEE) (S’98-M’01-SM’06-F’19) received the B.Tech. degree in electronics and communications engineering from the Indian Institute of Technology (IIT), Madras in 1996, and the M.S. and Ph.D. degrees in electrical engineering from the California Institute of Technology, Pasadena, USA, in 1997 and 2001, respectively. He is currently a Professor with the Department of Electrical Communication Engineering, Indian Institute of Science, Bengaluru. His research group works on the design, modeling, analysis, and optimization of current and next-generation wireless communication systems.
\end{IEEEbiography}

\begin{IEEEbiography}[{\includegraphics[width=1in,height=1.25in,clip,keepaspectratio]{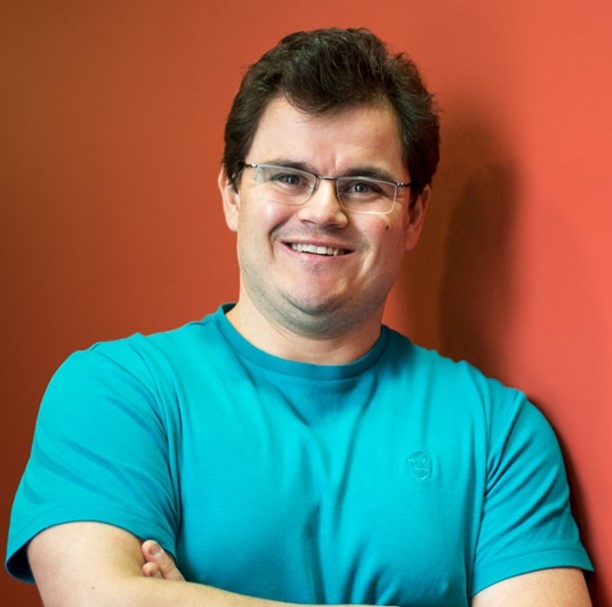}}]{Hirley Alves } (Member, IEEE) (S’11–M’15) received the B.Sc. and M.Sc. degrees from the Federal University of Technology-Paraná (UTFPR), Brazil, in 2010 and 2011, respectively in electrical engineering, and the dual D.Sc. Degree from the University of Oulu and UTFPR, in 2015. In 2017, he was an Adjunct Professor in machine-type wireless communications with the Centre for Wireless Communications (CWC), University of Oulu, Oulu, Finland. He is an Associate Professor at CWC and Leader of the Machine-type Wireless Communications Group and leads the Strategic Research Area 1 Massive Wireless Automation in the 6G Flagship. He is actively working on massive connectivity and ultra-reliable low latency communications for future wireless networks, 5GB and 6G, machine learning for wireless. He has been the organizer, chair, TPC, and tutorial lecturer for several renowned international conferences.
\end{IEEEbiography}

\begin{IEEEbiography}[{\includegraphics[width=1in,height=1.25in,clip,keepaspectratio]{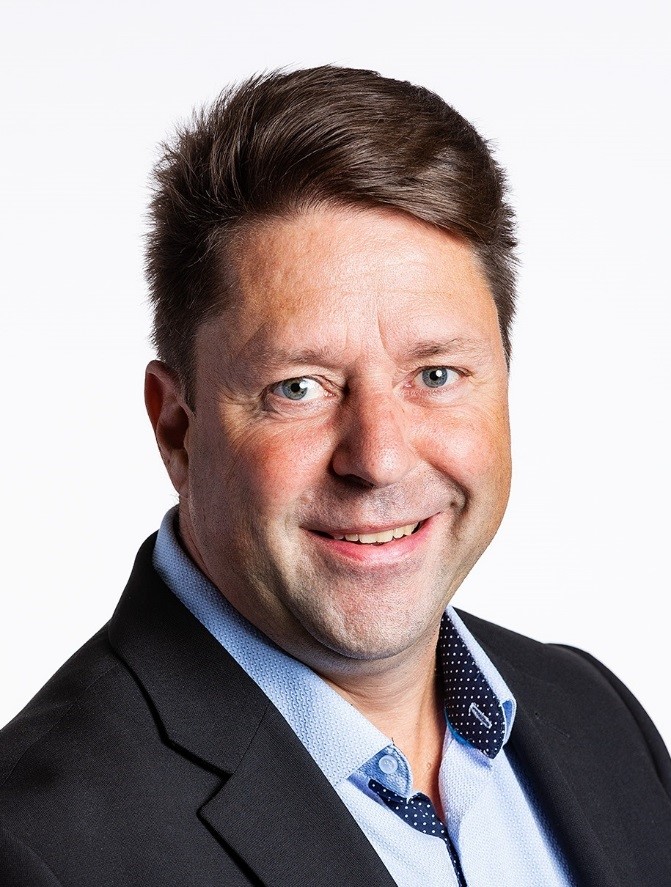}}]{Matti Latva-aho } (Fellow, IEEE) received the M.Sc., Lic.Tech. and Dr. Tech (Hons.) degrees in Electrical Engineering from the University of Oulu, Finland in 1992, 1996 and 1998, respectively. From 1992 to 1993, he was a Research Engineer at Nokia Mobile Phones, Oulu, Finland after which he joined the Centre for Wireless Communications (CWC) at the University of Oulu. Prof. Latva-aho was Director of CWC during the years 1998-2006 and Head of the Department for Communication Engineering until August 2014. Currently, he is a professor at the University of Oulu on wireless communications and Director for the National 6G Flagship Programme. He is also a Global Research Fellow with Tokyo University. His research interests are related to mobile broadband communication systems and currently, his group focuses on 6G systems research. Prof. Latva-aho has published over 500 conference or journal papers in the field of wireless communications. He received Nokia Foundation Award in 2015 for his achievements in mobile communications research.
\end{IEEEbiography}

\end{document}